\documentclass[11pt]{article} 
\usepackage{amscd,amsmath, amssymb, fancyhdr, epsfig,color}
\usepackage{hyperref}

\usepackage[small,nohug,heads=vee]{diagrams}
\diagramstyle[labelstyle=\scriptstyle]

\numberwithin{equation}{section}

\usepackage{letltxmacro}
\makeatletter
\let\oldr@@t\r@@t
\def\r@@t#1#2{%
\setbox0=\hbox{$\oldr@@t#1{#2\,}$}\dimen0=\ht0
\advance\dimen0-0.2\ht0
\setbox2=\hbox{\vrule height\ht0 depth -\dimen0}%
{\box0\lower0.4pt\box2}}
\LetLtxMacro{\oldsqrt}{\sqrt}
\renewcommand*{\sqrt}[2][\ ]{\oldsqrt[#1]{#2} }
\makeatother

\newcommand{\Pic}{\mathrm{Pic}}
\DeclareMathOperator{\Hilb}{Hilb}
\DeclareMathOperator{\Sym}{Sym}
\DeclareMathOperator{\Vol}{Vol}

\newcounter{Mycounter}[section]
\newcounter{lemma}[section]
\newcounter{claim}[section]
\newcounter{sublemma}[section]
\newcounter{corollary}[section]
\newcounter{remark}[section]
\newcounter{theorem}[section]
\renewcommand{\thetheorem}{{Theorem \thesection.\arabic{theorem}}}
\newcommand{\theorem}{%
     \setcounter{theorem}{\value{Mycounter}}
     \refstepcounter{theorem}
     \stepcounter{Mycounter}
     {\bf \thetheorem:\ }}

\newcounter{conjecture}[section]
\newcounter{proposition}[section]
\newcounter{definition}[section]
\renewcommand{\thedefinition}
       {{Definition~\thesection.\arabic{definition}}}
\newcommand{\definition}{%
     \setcounter{definition}{\value{Mycounter}}
     \refstepcounter{definition}
     \stepcounter{Mycounter}
     {\bf \thedefinition:\ }}
     
\newenvironment{proof}{
{\bf Proof:}}{$\blacksquare$}

\begin{document}
%%%%%%%%%%%%%%%%%%%%%%%%%%%%%%%%%%%%%%%%%%%%%%%%%%%%%%%%%%%%
\begin{center}
{\LARGE\bf
Absolutely trianalytic tori in the generalized Kummer variety\\[3mm]
}
%%%%%%%%%%%%%%%%%%%%%%%%%%%%%%%%%%%%%%%%%%%%%%%%%%%%%%%%%%%%

Nikon Kurnosov\footnote{Nikon Kurnosov is partially supported by grant MK-1297.2014.1. The article was prepared within the framework of a subsidy granted to the HSE by the Government of the Russian Federation for the implementation of the Global Competitiveness Program}

\end{center}

%%%%%%%%%%%%%%%%%%%%%%%%%%%%%%%%%%%%%%%%%%%%%%%%
{\small \hspace{0.10\linewidth}
\begin{minipage}[t]{0.85\linewidth}
{\bf Abstract} \\
We prove that a generic complex deformation of a generalized Kummer variety contains no complex analytic tori.
\end{minipage}
}
%%%%%%%%%%%%%%%%%%%%%%%%%%%%%%%%%%%%%%%%%%%%%%%%

\tableofcontents

%%%%%%%%%%%%%%%%%%%%%%%%%%%%%%%%%%%%%%%%%%%%%%%%

\section{Introduction}

A Riemannian manifold is called \textit{hyperk\"ahler} if it admits a triple of complex structures
$I,J,K$ satisfying quaternionic relations and K\"ahler with
respect to $g$. 

\hfill

\definition \label{Def_holom_symplectic} 
A manifold $M$ is called \textit{holomorphically symplectic} if it is a complex manifold with a closed holomorphic 2-form $\Omega$ over $M$ such that 
$\Omega^n = \Omega \wedge \Omega \wedge ... \wedge \Omega$ is
a nowhere degenerate section of a canonical class of $M$, where $2n = \dim_\mathbb{C}(M)$.

\hfill

A hyperk\"ahler manifold is always holomorphically symplectic. By the Yau's Theorem \cite{Y}, a hyperk\"ahler structure exists on a compact complex manifold 
if and only if it is K\"ahler and holomorphically symplectic.

Given any triple $a, b ,c \in \mathbb{R}$, $a^2+b^2+c^2=1$,
the operator $L:=aI+bJ+cK$ satisfies $L^2=-1$ and defines 
a K\"ahler structure on $(M,g)$. Such a complex structure is called 
\textit{induced by the hyperk\"ahler structure}. Complex subvarieties
of such $(M,L)$ for generic $(a,b,c)$ were studied in \cite{V1}, \cite{V2}.

\definition \label{_IHS_Definition_}
A compact hyperk\"ahler manifold $M$ is called 
\textit{simple}, or \textit{irreducible holomorphically symplectic (IHS)} or \textit{of maximal holonomy} if $\pi_1(M)=0$, $H^{2,0}(M)=\mathbb{C}$.

\hfill

\theorem \label{_Bogo_decompo_Theorem_}
(Bogomolov's decomposition, \cite{B})
Any  hyperk\"ahler manifold admits a finite covering
which is a product of a torus and several 
 hyperk\"ahler manifolds of maximal holonomy.

\hfill

\definition\label{_twistor_family_Definition_}
Let $M$ be a hyperk\"ahler manifold, and all induced complex structures $L:=aI+bJ+cK$, where
$a, b ,c \in \mathbb{R}$, $a^2+b^2+c^2=1$ fit together into a family over $\mathbb{C}P^1$ called \textit{the twistor family} of complex structures.

\hfill

\definition
A closed subset $Z$ of a 
hyperk\"ahler manifold $M$ is called \textit{trianalytic}
if it is complex analytic with respect to complex structures $I, J, K$.

\hfill

\definition \label{_abs_triana_Definition_}
Let $(M,I,J,K)$ be a compact, holomorphically symplectic,
K\"ahler manifold, and $Z\subset (M,I)$ a complex subvariety
which is trianalytic with respect to any hyperk\"ahler structure
compatible with $I$. Then $Z$ is called \textit{absolutely trianalytic}.

\hfill

Whenever $L$ is a generic
element of a twistor family, all subvarieties
of $(M,L)$ are trianalytic (see \ref{generic}).

Absolutely trianalytic subvarieties were studied 
in \cite{V3}
, where it was shown that
a general deformation of a Hilbert scheme of a $K3$ surface
has no complex (or, equivalently, no absolutely trianalytic (\cite[Theorem 8.5.]{V3}))
subvarieties. However, there is an absolutely trianalytic subvariety in a generalized Kummer variety. Recently, Soldatenkov and Verbitsky \cite{SV} have shown that there are no absolutely trianalytic tori in the 6- and 10-dimensional O'Grady examples. Non-existence of absolutely trianalytic subvarieties
of known type in 10-dimensional O'Grady manifold $M$ follows from \cite[Corollary 3.17]{SV}. Non-existence of absolutely trianalytic tori in a 6-dimensional O'Grady manifold follows from representation theory of Clifford algebras \cite{SV}.
Automorphisms of hyperk\"ahler manifolds acting trivially on the second cohomology group can be characterized as those which have absolutely trianalytic graph in $M \times M$ for any given hyperk\"ahler structure on $M$. The group of such automorphisms is finite \cite{H1}. It has been studied for the generalized Kummer surfaces by Oguiso (\cite{Og}), Boissiere, Nieper-Wisskirchen and Sarti (\cite{BNS}), and for O'Grady examples by Mongardi and Wandel (\cite{MW}).

In the present paper we study absolutely trianalytic tori in the generalized Kummer variety. In Section \ref{Sec2} we recall all known examples of hyperk\"ahler manifolds. In Section \ref{Sec3} we study general trianalytic subvarieties. In Section \ref{Sec4} we show non-existence of absolutely trianalytic tori in the generalized Kummer variety and prove the following

\hfill

\begin{theorem}\label{non-exist}
Let $K_n(T)$ be a generalized Kummer variety, and $Z \subset K_n(T)$ be an absolutely trianalytic subvariety of $K_n(T)$. Then $Z$ is not a torus.
\end{theorem}

\section{Preliminaries} \label{Sec2}

Let $(M, I, J, K)$ be a hyperk\"ahler manifold, and let $\omega_I, \omega_J, \omega_K$ be the corresponding K\"ahler forms.

\hfill
 
A simple algebraic calculation \cite{Bes} shows that the following form
\begin{equation}
\Omega=\omega_J+\sqrt{-1}\omega_K
\end{equation}
is of type $(2,0)$. Since it is closed this form is also holomorphic and moreover nowhere degenerate, as another linear algebraic argument shows.
It is called \textit{the canonical holomorphic symplectic form
of a manifold $M$}. Thus, the underlying complex manifold
$(M,L)$ is holomorphically symplectic for each hyperk\"ahler manifold $M$
and an induced complex structure $L$. The converse is also true:

\hfill

\theorem \label{symplectic=>hyperkahler}
(\cite{Bea}, \cite[Chapter 11]{Bes})
Let $M$ be a compact holomorphically
symplectic K\"ahler manifold with the 
holomorphic symplectic form
$\Omega$, a K\"ahler class 
$[\omega]\in H^{1,1}(M)$ and a complex structure $I$. 
Let $n=\dim_\mathbb{C} M$. Assume that
$\int_M \omega^n = \int_M ($Re$ \Omega)^n$.
Then there is a unique hyperk\"ahler 
structure $(I,J,K,(\cdot,\cdot))$
over $M$ such that the cohomology class of the symplectic form
$\omega_I=(\cdot,I\cdot)$ is equal to $[\omega]$ and the
canonical symplectic form $\omega_J+\sqrt{-1}\omega_K$ is
equal to $\Omega$.

\hfill

Two-dimensional irreducible holomorphic symplectic manifolds are $K3$ surfaces. In higher dimensions there are
only few examples known. Here is the list of known examples, where compact manifolds of
the same deformation type are not distinguished.

\begin{itemize}
\item[(0)] \textit{$K3$ surface.}

\item[(i)] \textit{The Hilbert scheme of $n$ points of $K3$.} If $X$ is a $K3$ surface then the Hilbert scheme $\Hilb^n(X)$ is an irreducible
holomorphic symplectic manifold \cite{Bea}. Its dimension is $2n$ and for $n > 1$ its
second Betti number is equal to 23. Details of construction of the Hilbert scheme can be
found, for example \cite{Bea}. Namely, let $X$ be a $K3$ surface. Take
the symmetric product $X^{\left( r \right)} = X^r /\mathfrak{S}_r$ which
parametrizes subsets of $r$ points in a $K3$ surface $X$, counted with multiplicities;
it is smooth on the open subset $X_0$ consisting of subsets with $r$ distinct points,
but singular otherwise. We blow up singular locus and obtain a smooth compact manifold. This is 
the Hilbert scheme $X^{\left[ r \right]}$. The natural map $X^{\left[ r \right]} \rightarrow X^{\left( r
\right)}$ is an isomorphism above $X_0$, and it resolves the singularities of $X^{\left( r \right)}$. Alternatively, the Hilbert scheme parametrizes all $0$-dimensional subschemes of the length $n$.

Let us describe the simplest case $\Hilb^2(X)$ explicitly. For any surface $X$ the
Hilbert scheme $\Hilb^2(X)$ is the blow-up $\Hilb^2 (X) \rightarrow S^2
(X)$ of the diagonal 
\[\Delta =\{\{x, x\} \left| \right. x \in X\} \subset S^2
(X) =\{\{x, y\} \left| \right. x, y \in X\}.\]
Equivalently, $\Hilb^2(X)$ is the
$\mathbb{Z}/ 2\mathbb{Z}$-quotient of the blow-up of the diagonal in $X
\times X$. Since for a $K3$ surface there exists only one $\mathbb{Z}/
2\mathbb{Z}$-invariant two-form on $X \times X$, the holomorphic symplectic
structure on $\Hilb^2(X)$ is unique.

\item[(ii)] \textit{The generalized Kummer variety.} If $T$ is a complex torus of dimension two, then the generalized Kummer
variety K$_n(T)$ is an irreducible holomorphic symplectic manifold \cite{Bea}. Its
dimension is $2n$ and for $n > 2$ its second Betti number is 7. Note that the Hilbert
scheme $T^{\left[ n \right]}$ of a two-dimensional torus has the same properties
as $K3^{\left[ r \right]}$, but it is not simply connected. The commutative group structure on the torus $T$ defines a summation map
\[s \left( t_1,
\ldots, t_n \right) = t_1 + \ldots + t_{n+1},\]
\[\Sigma: T^{n+1} \rightarrow T,\]
 which induces a morphism $\Sigma: T^{[n+1]} \rightarrow T$. It is easy to see
that $\Sigma$ coincides with the Albanese map. The generalized Kummer variety K$_n(T)$ associated to the
torus $T$ is the preimage $\Sigma^{-1}(0) \subset T^{[n+1]}$ of the zero $0 \in T$. It is a hyperk\"ahler
manifold of dimension $2n$.

\item[(iii)] \textit{O'Grady's 10-dimensional example \cite{O1}.} Let again $X$ be a $K3$ surface, and
M the moduli space of stable rank 2 vector bundles on $S$, with Chern classes
$c_1 = 0, c_2 = 4$. It admits a natural compactification $M$ obtained by adding classes of semi-stable torsion free sheaves. It is singular
along the boundary, but O'Grady \cite{O1} constructs a desingularization of $M$ which
is a new hyperk\"ahler manifold, of dimension 10. Its second Betti number is 24
\cite{R}. Originally, it was proved that it is at least 24 \cite{O1}.

\item[(iv)] \textit{O'Grady's 6-dimensional example \cite{O2}.} A similar construction can be
done starting from rank 2 bundles with $c_1 = 0, c_2 = 2$ on a 2-dimensional
complex torus, this gives new hyperk\"ahler manifold of dimension 6 as in (iii). Its second Betti number is 8.

\end{itemize}

Thus we have two series, (i) and (ii), and two sporadic examples, (iii) and
(iv). All of them have different second Betti numbers. It has been proved (\cite{KLS}, Theorem B) that up to a deformation the moduli spaces for all sets of numerical parameters give $\Hilb^n(K3)$, O'Grady examples, or do not admit a smooth symplectic resolution of singularities. 

\section{Subvarities in hyperk\"ahler manifolds} \label{Sec3}

\subsection{Trianalytic subvarities}

\definition
A closed subset $Z$ of a 
hyperk\"ahler manifold $M$ is called \textit{trianalytic}
if it is complex analytic with respect to complex structures $I, J, K$.

\hfill

\begin{theorem} (equivalent to the \ref{symplectic=>hyperkahler}) \\Let $M$ be a hyperk\"ahler manifold. Then there exists a unique hyperk\"ahler metric in a given K\"ahler class.
\end{theorem}

\hfill

\definition \label{_abs_triana_Definition_}
Let $(M,I,J,K)$ be a compact, holomorphically symplectic,
K\"ahler manifold, and $Z\subset (M,I)$ a complex subvariety
which is trianalytic with respect to any hyperk\"ahler structure
compatible with $I$. Then $Z$ is called \textit{absolutely trianalytic}.

\hfill

\begin{theorem} \cite{SV}\label{defM-M'}
For any hyperk\"ahler manifolds $M, M'$ in the same deformation class there is a diffeomorphism which sends absolutely trianalytic subvarieties to absolutely trianalytic.
\end{theorem}

\hfill

\begin{definition}\label{general}
A hyperk\"ahler manifold is called \textit{general} if all its subvarities are absolutely trianalytic.
\end{definition}

\hfill

%\begin{definition}\label{generic}
%Let $M$ be a compact hyperk\"ahler manifold, and $I$ an induced complex structure. We say that $I$ is \textit{of general type} or \textit{generic} with respect to the hyperk\"ahler structure on $M$, if all elements of the group
%\[\bigoplus_p H^{p,p}(M) \cap H^{2p}(M, \mathbb{Z}) \subset H^*(M)\]
%are $SU(2)$-invariant.
%\end{definition}

%\hfill

\begin{remark} General deformation of a hyperk\"ahler manifold is general in the sense of \ref{general} (\cite[Proposition 2.14]{KV-book}).
\end{remark}

\hfill

\begin{theorem} \label{generic} (\cite{V2})
Let $M$ be a hyperk\"ahler manifold, $S$ its twistor family (see \ref{_twistor_family_Definition_}).
Then there exists a countable subset $S_1\subset S$, such that
for any complex structure $L\in S \backslash S_1$, all
compact complex subvarieties of $(M,L)$ are trianalytic.
\end{theorem}

\subsection{Trianalytic subvarities in the Hilbert scheme and O'Grady examples}

Here we survey the known results about trianalytic and absolutely trianalytic subvarieties.

It was shown by Verbitsky that a general deformation of a Hilbert scheme of a $K3$ surface has no complex subvarieties \cite{V3}. The same theorem was also claimed (Kaledin, Verbitsky) in the case of generalized Kummer varieties \cite{KV}. However, later (\cite{KV1}) they found that there are counterexamples in the latter case due to involution $\nu: t \rightarrow -t$ of a torus. This involution is extended to an involution of the Hilbert scheme $T^{[n+1]}$, and since it commutes with the Albanese map $T^{[n+1]} \longrightarrow T$, the map $\nu$ preserves $K_n(T)$. Moreover, $\nu$ sends the K\"ahler class to itself. Hence, the involution $\nu$ preserves the hyperk\"ahler structure on $K_n(T)$. For odd $n = 2m - 1$ the map $\nu$ fixes the $2m$-tuple
\[(x_1, - x_1, x_2, -x_2, ..., x_m, -x_m) \in T^{(n+1)}\]
When $x_i , -x_i$ are pairwise distinct, they give a point of the Hilbert scheme fixed by $\nu$. Consider the closure $X$ of the set of such points. It is one of components of fixed point set of involution map $\nu$. The submanifold $X$ is birationally equivalent to the Hilbert scheme of a $K3$ surface.

\hfill

\begin{corollary} The variety $\Sym^2(T)$ contains a Kummer $K3$ surface.
\end{corollary}

\hfill

Non-existence of absolutely trianalytic subvarieties
in the Hilbert scheme $\Hilb(K3)$ of $K3$ was used in the book \cite{KV-book} to prove compactness of 
deformation spaces of certain stable holomorphic bundles
on $M$.

\hfill

\begin{theorem}(\cite[Theorem 6.2.]{V5}) 
Let $M$ be a hyperk\"ahler manifold, $Z\subset M$ a trianalytic
subvariety, and
$I$ an induced complex structure.
Consider the normalization \[ \widetilde{(Z, I)} \rightarrow (Z,I)\] 
of $(Z,I)$. Then $\widetilde{(Z, I)}$ is smooth, and
the map $\widetilde{(Z, I)} \rightarrow  M$ is an immersion,
inducing a hyperk\"ahler structure on $\widetilde{(Z, I)}$.
\end{theorem}

\hfill

This gives that any trianalytic subvariety $Z \rightarrow M$
has a smooth hyperk\"ahler normalization $\widetilde{Z}$ immersed to $M$;
this immersion is generically bijective onto its image. Therefore, we can replace any trianalytic
cycle by an immersed hyperk\"ahler manifold. In this paper, we will consider absolutely trianalytic varieties whose normalization is the torus.

\hfill

\begin{theorem} \cite{SV} Let $M$ be a hyperk\"ahler manifold,
$Z\subset M$ an absolutely trianalytic subvariety, and
$\widetilde{Z} \rightarrow M$ its normalization such that $\widetilde{Z} = T\times\prod_i K_i$, where $K_i$ are IHS of maximal holonomy. 
Then $b_2(T)\geq b_2(M)$ and $b_2(K_i)\geq b_2(M)$.
\end{theorem}

\hfill

\begin{theorem} \label{ogrady} \cite{SV} Let $M$ be a hyperk\"ahler manifold of maximal holonomy, 
$T$ a hyperk\"ahler torus, and $T\rightarrow M$ a hyperk\"ahler
immersion with absolutely trianalytic image. Then 
\[\dim_{\mathbb{C}}(T) \geqslant 2^{\frac{b_2(X)-1)}{2}}.\]
\end{theorem}

\section{Tori in the generalized Kummer varieties} \label{Sec4}

In this section we prove the following

\textbf{Main Theorem (\ref{non-exist}:)}
Let $K_n(T)$ be a generalized Kummer variety, and let $Z \subset K_n(T)$ be an absolutely trianalytic subvariety. Then $Z$ is not a torus.

%\hfill

\subsection{Flat tori in $T^n$} \label{subsecDef}

Let $Z \subset K_n(T)$ be an absolutely trianalytic subvariety whose normalization is a torus.

\hfill

\begin{remark} Since $K_n(T)$ is embedded in the Hilbert scheme of a torus $T^{[n]}$, we can consider $Z$ as an absolutely trianalytic submanifold in $T^{[n]}$.
\end{remark}

\hfill

Consider the diagram:

\begin{displaymath}
\begin{diagram}[labelstyle=\scriptscriptstyle]
Z& \rTo & \pi(Z) & \lTo & \tau^{-1}(\pi(Z))\\
\dTo & \rdTo & \dTo &  & \dTo \\
T^{[n]} & \rTo^{\pi} & T^{(n)} &  \lTo^{\tau} & T^n,
\end{diagram}
\qquad \eqno {(4.1)} \label{diag4.1}
\end{displaymath}

where $T^{[n]}$ is the Hilbert scheme of a torus, $T^{(n)}$ is the symmetric power of a torus, the map $\pi$ is the Hilbert-Chow map, $\tau$ is the quotient map $T^n \rightarrow T^{(n)}$, and the square is Cartesian.

\hfill

\begin{remark} By \cite[Theorem 7.7]{EV} one can choose $T$ in the same deformation class that is general in the sense of \ref{general}. Recall that by \ref{defM-M'} trianalytic subvarieties have the following property: for each $M, M'$ in the same deformation class there exist diffeomorphism $M \to M'$ which sends absolutely trianalytic subvarities to absolutely trianalytic ones. Therefore it is sufficient to prove our proposition for a general torus $T$.
%The projection $Z \rightarrow \pi(Z) \subset T^{(n)}$ will be denoted by $\pi$.
\end{remark}

\hfill

\begin{proposition}\label{tau-pi-Z}
Let $Z \subset T^{[n]}$ be an absolutely trianalytic torus. Then each irreducible component of $\tau^{-1}(\pi(Z))$ is a general torus in $T^n$, $\Pic(\tau^{-1}(\pi(Z))) = 0$, and the map $\tau: \tau^{-1}(\pi(Z)) \rightarrow \pi(Z)$ is finite.
\end{proposition}

\begin{proof}
Recall that $T$ is a general torus. Therefore, all subvarieties in $T^n$ are absolutely trianalytic, hence $\tau^{-1}(\pi(Z))$ is totally geodesic, and hence is flat. Therefore each irreducible component of $\tau^{-1}(\pi(Z))$ is a subtorus in $T^{n}$. Since $T$ is general, any subtorus is general, and in particular the Picard group of $\tau^{-1}(\pi(Z))$ is zero. The finiteness of the map $\tau$ is clear.
\end{proof}

\hfill

\begin{proposition}\label{finite-sym}
Let $Z \subset T^{[n]}$ be an absolutely trianalytic torus. Then the map $\pi: Z \rightarrow \pi(Z)$ is generically finite.
\end{proposition}

\begin{proof}
There exist a canonical stratification on every symplectic singularity and this stratification coincide with stratification by diagonals on $T^{(n)}$ (\cite[Proposition 3.1.]{K}). Moreover, strata carry symplectic forms, and for $T^{(n)}$ these forms are induced by the transposition-equivariant symplectic form on $T^n$. The restriction of the symplectic form to the smooth locus of the preimage $\pi^{-1}(V)$ for an arbitrary stratum $V$ in $T^{(n)}$ is the pullback of a symplectic form on this stratum in $T^{(n)}$ (\cite[Lemma 2.9.]{K}). Then a dense open subset $U$ in $Z \subset T^{[n]}$ is projected into the open part of some stratum. Thus, the restriction of the symplectic form to $U$ is a pullback of the symplectic form on the stratum. If $Z$ is projected to $\pi(Z)$ with a positive-dimension general fibres, then the form cannot be non-degenerate, so that $Z$ is not symplectic. This gives a contradiction.
\end{proof}

%эквиваоиантна поэтому спускается.
\hfill

\begin{proposition}\label{flat-tori}
Let $Z \subset T^{[n]}$ be an absolutely trianalytic torus. Consider the diagram

\begin{displaymath}
\begin{diagram}[labelstyle=\scriptscriptstyle]
& &\tilde{Z}\\
& \ldTo &  & \rdTo & \\
Z &  &  &  & \tau^{-1}(\pi(Z))\\
& \rdTo &  & \ldTo & \\
& &\pi(Z)\\
\end{diagram}
\qquad \eqno {(4.2)} \label{diag4.2}
\end{displaymath}

where $\tilde{Z}$ is the fibered  product of $Z$ and $\tau^{-1}(\pi(Z))$. Then $Z$ and any component of $\tau^{-1}(\pi(Z))$ are isogenic tori.
\end{proposition}

\begin{proof}

The fibered product $\tilde{Z}$ is a subvariety in the product $Z\times\tau^{-1}(\pi(Z)$ of general tori. Then $\tilde{Z}$ is trianalytic, and therefore flat. Maps of flat tori $\tilde{Z}$ to $Z$ and $\tau^{-1}(\pi(Z))$ are generically finite (finiteness of a map $\tilde{Z} \rightarrow \tau^{-1}(\pi(Z))$ follows from \ref{finite-sym}), hence these tori are isogenic.
\end{proof}

\hfill

Fix an irreducible component of $\tau^{-1}(\pi(Z))$ and denote it by $Z'$.

Denote the generic degree of the map $Z \rightarrow \pi(Z)$ by $d$ and the degree of $\tau^{-1}(\pi(Z)) \rightarrow \pi(Z)$ by $\tilde{d}$.

\hfill

\begin{remark} \label{Pic} It follows from \ref{tau-pi-Z} and \ref{flat-tori} that $\Pic(Z) = 0$
\end{remark}

\subsection{Non-existence of trianalytic tori in the generalized Kummer variety}

In this section we prove the Main \ref{non-exist}.

\hfill

\begin{definition}
The holomorphic symplectic volume of a holomorphic symplectic manifold $(M, \Omega)$ is $\Vol^s_M := \int_M\Omega^{\frac{1}{2}\dim M}\wedge\overline{\Omega}^{\frac{1}{2}\dim M}$.
\end{definition}

\hfill

\begin{definition}
The K\"ahler volume of a K\"ahler manifold $(M, I, \omega)$ is $\frac{1}{2^{2n}(2n)!}\int_M\omega^{2n}$, where dim$_\mathbb{R}(M) = 2n$.
\end{definition}

\hfill

\begin{remark} For hyperk\"ahler manifolds (\cite[Theorem 5.3.]{GV}) the K\"ahler volume is equal to the holomorphic symplectic one (the \textit{hyperk\"ahler condition}).
\end{remark}

\hfill

\begin{theorem}
Let $K_n(T)$ be a generalized Kummer variety, and let $Z \subset K_n(T)$ be an absolutely trianalytic subvariety of $K_n(T)$. Then $Z$ is not a torus.
\end{theorem}

\begin{proof}

First, let us remark that any complex
structure of K\"ahler type on a flat torus $T$ defines a complex structure of K\"ahler type on
$T^{[n]}$. Consider the standard map from $H^2(T^{(n)},\mathbb{C})\oplus \mathbb{C}[E]$ to $H^2 (T^{[n]},\mathbb{C})$, where $E$ is the exceptional divisor of the blow-up $T^{(n)}$ to $T^{[n]}$. Recall that the cohomology class $[E]$ is of type $(1,1)$.

Fix a hyperk\"ahler structure ($I, J, K$) on $T$, and let $\Omega$ be the corresponding holomorphic symplectic form on $T^n$. Hyperk\"ahler triple on $T^n$ is given by three forms $\omega'_I$, $\omega'_J$, and $\omega'_K$ of the same K\"ahler volume. Denote by $[\omega'_I]$, $[\omega'_J]$, and $[\omega'_K]$ their cohomology classes. Since the symplectic form on $T^n$ is transposition-equivariant, the corresponding cohomology classes on $T^{(n)}$ denoted by $[\omega_I]$, $[\omega_J]$, and $[\omega_K]$  are such that
\[ \tau^*[\omega_I] = [\omega'_I], \tau^*[\omega_J] = [\omega'_J], \tau^*[\omega_K] = [\omega'_K]. \]
By \ref{finite-sym} on each open stratum of $T^{(n)}$ there exist a symplectic form, and the cohomology class of this form is the restriction of $[\omega_J] + i[\omega_K]$ to the stratum.
 %Since $\Omega$ is invariant under permutations, we can consider $\Omega$ as a holomorphic symplectic form on the orbifold $T^{(n)}$. Therefore, we could make computations for 2-forms whose are well-defined.

It is well-known (see e.g. \cite[Lemma 3.4]{OVV}) that there exist a K\"ahler metric with K\"ahler class $[\omega_{T^{[n]}}] = [\pi^*\omega_{T^{(n)}}] - \epsilon [E]$, where $E$ is the exceptional divisor and $0 < \epsilon < 1$.

Recall that the symplectic volume does not change under the blow-up. By \ref{symplectic=>hyperkahler} there exists some constant $\lambda$ and a hyperk\"ahler structure on $T^{[n]}$ such that $[\tilde{\omega}_I] := \lambda [\pi^*\omega_I] - \lambda \epsilon [E]$, $[\tilde{\omega}_J]$, and $[\tilde{\omega}_K]$ have the same K\"ahler volume. After blowing up the symmetric power $T^{(n)}$ to the Hilbert scheme of points $T^{[n]}$ pullbacks of $[\omega_J]$ and $[\omega_K]$ are $[\tilde{\omega}_J]$ and $[\tilde{\omega}_K]$:
\[ \pi^*[\omega_J] = [\tilde{\omega}_J], \qquad \pi^*[\omega_K] = [\tilde{\omega}_K].\]

Note that $\pi^*(\omega)\cup[E] = 0$, and let $\mu = [E]^{2n}$.

Then,

\[\Vol_{\omega_I} = \int_{T^{[n]}} \tilde{\omega}_I^{2n} = (\lambda)^{2n} \cdot \Vol_{\omega_I} -\lambda \epsilon^{2n}\cdot\mu\]

The constant $\lambda$ can be determined from the equation above, and we have $\lambda >1$.

\hfill

%Whenever a manifold is hyperk\"ahler, its symplectic volume is equal to the K\"ahler one. 
%It follows from \ref{flat-tori} that the map from $Z$ to each component of $\tau^{-1}(\pi(Z)))$ is an isogeny. Thus, volumes of $Z$ and $\tau^{-1}(\pi(Z)))$ differs by the multiplication by $d$, where $d$ is the degree of the ramified covering $\tau^{-1}(\pi(Z))) \rightarrow Z$. The K\"ahler volume is determined from the hyperk\"ahler condition.
Recall that by \ref{finite-sym}, and \ref{tau-pi-Z}, the map from $Z$ to $\pi(Z)$ and the map from $Z'$ to $\pi(Z)$ are generically finite. Thus, symplectic volumes of $Z$ and $Z'$ differ by the multiplication by $\frac{\tilde{d}}{d}$
\[ \frac{\tilde{d}}{d} \cdot \Vol^s_Z = \Vol^s_{Z'},\]
where $d$ is the degree of the map $Z \rightarrow \pi(Z)$ and $\tilde{d}$ is the generic degree of the map $\tau^{-1}(\pi(Z)) \rightarrow \pi(Z)$. The K\"ahler volumes of $Z$ and $Z'$ are determined from the hyperk\"ahler condition.

%Let us calculate the K\"ahler volumes of $Z$ in $T^{[n]}$ and $\tau^{-1}(\pi(Z))$ in $T^n$.
Since $Z'$ is absolutely trianalytic in $T^{n}$, its volume with respect to $\omega'_J$ and $\omega'_K$ is equal to the volume with respect to $\omega'_I$. However, $Z$ is also absolutely trianalytic in $T^{[n]}$, hence this volume is also equal to the volume with respect to $\lambda [\pi^*\omega_I] - \lambda \epsilon [E]$.

Note that $[Z]\cup[E] = 0$. Indeed, consider line bundle $\mathcal{O}(E)$ restricted to $Z$. Since $Z$ has zero Picard group (\ref{Pic}), then there are no non-trivial line bundles over $Z$.\\

Hence, we have from the hyperk\"ahler condition

\[ 1 = \frac{\int_Z (\lambda[\pi^*\omega_I] - \lambda \epsilon [E])^k}{\int_Z (\tilde{\omega}_J)^k} = \frac{\lambda^k \int_{\pi(Z)} [\omega_I]^k}{\int_{\pi(Z)} [\omega_J]^k} = \frac{\lambda^k \int_{Z'} (\omega'_{I})^k}{\int_{Z'} (\omega'_{J})^k} = \lambda^k.
\]

On the other hand $\lambda > 1$, that gives a contradiction.
\end{proof}

\hfill

\begin{remark} In the general case of absolutely trianalitic subvarieties of generalized Kummer manifold the proof above does not work. Indeed, generally $Z$ is not isogenic to irreducible components of $\tau^{-1}(\pi(Z))$, and Picard group of $Z$ could be non-zero.
\end{remark}

\hfill

%\corollary Let $Z$ be an absolutely trianalytic subvariety of the generalized Kummer $K_n(T)$. Then $\Pic(Z) \neq 0$.

%\hfill

{\bf Acknowledgement:}
Author is grateful to his supervisor Misha Verbitsky for constant interest and fruitful
discussions, to Andrey Soldatenkov for comments, to Dmitry Kaledin for interesting discussions and suggestions, to Constantin Shramov, who read the draft of the paper and gives useful comments, and to the reviewer for important corrections.

\hfill

{\small
\noindent {\sc Nikon Kurnosov\\
Laboratory of Algebraic Geometry and its applications,\\
National Research University Higher School of Economics\\
7 Vavilova Str., Moscow, Russia, 117312;\\
Independent University of Moscow,\\
11 Bol.Vlas'evskiy per., Moscow, Russia, 119002}\\
\it  nikon.kurnosov@gmail.com 
}

\end{document}